\documentclass[11pt]{article}
\usepackage{amsthm, amsmath,amssymb,tabularx,amsthm}
\usepackage{longtable}
\usepackage[T1]{fontenc}
\usepackage[latin1]{inputenc}
\newtheorem{theorem}{Theorem}[section]
\newtheorem{corollary}[theorem]{Corollary}
\newtheorem{lemma}[theorem]{Lemma}
\usepackage{url} 
\usepackage{icomma} 
\usepackage{calc}
\usepackage[hang]{caption}
\usepackage{varioref}
\usepackage{booktabs,longtable}
\numberwithin{equation}{section}
\usepackage{algpseudocode}
\usepackage{MnSymbol}
\usepackage{multirow}
\algrenewcommand\algorithmicindent{0.7em}

\begin{document}

	\title{New results on maximal partial line spreads in $\mathrm{PG}(5, q)$}
	
	\author{Maurizio Iurlo}
	
	\date{\,}
	
	\maketitle

	\begin{abstract}
		In this work, we prove the existence of maximal partial line spreads in
		$\mathrm{PG}(5, q)$ of size $q^{3}+q^{2}+kq+1$, with $1\leq k\leq
		\left\lfloor \frac{q^{3}-q^{2}}{q+1}\right\rfloor $,
		$k$ an integer. Moreover, by a computer search, we do this for larger values of $k$, for $q\leq7$.
		Again by a computer search, we find the sizes for the
		largest maximal partial line spreads 
		and many new results for $q\leq5$.
			
	\end{abstract}
	
	
	\vspace{0,2cm}\noindent \textbf{AMS Classification:}\,\,\,51E14
	
	\section{Introduction}
	A \emph{partial line spread} $\mathcal{F}$ in $\mathrm{PG}(N, q)$,
	the projective space of dimensions $N$  over the Galois
	field $\mathrm{GF}(q)$ of order $q$,
	is a set of pairwise
	skew lines. We say that $\mathcal{F}$ is \emph{maximal} if it cannot be extended to a larger
	partial line spread. A \emph{line spread} in $\mathrm{PG}(N, q)$, $N$ odd, is a set of pairwise skew lines covering the space.

Here, we consider  a hyperplane $\mathcal{H}$ of $\mathrm{PG}(5, q)$ and 
	a  partial line spread  $\mathcal{F}$ of $\mathcal{H}$, of size 
	$q^{3}+1$, which is the largest size of MPS in $\mathrm{PG}(4, q)$. 
	So, in $\mathcal{H}$ there are exactly $q^{2}$\emph{ holes},
	that is points not over the lines of $\mathcal{F}$. 
	
	We start
	by adding to $\mathcal{F}$ a set of $q^2$ pairwise skew lines, not of $\mathcal{H}$, covering the holes of $\mathcal{H}$, obtaining a maximal partial spread in $\mathrm{PG}(5, q)$ of size $q^3+q^2+1$. We proceed
	by depriving the obtained maximal partial spread in $\mathrm{PG}(5, q)$ of some lines of $\mathcal{H}$ and adding $q+1$ pairwise skew
	lines not of $\mathcal{H}$ for each removed line. We do this through a
	theoretical way for every value of $q$, and by a computer search for
	$q\leq 7$. More precisely we prove that for every
	$q$ there are maximal partial line spreads of the sizes $q^{3}+q^2+kq+1$, for any $1 \leq k \leq 
	\left\lfloor \frac{q^{3}-q^{2}}{q+1}\right\rfloor$, $k$ an integer, while by a computer search we get cardinalities $q^{3}+q^{2}+kq+1$,
	for larger values of $k$, for $q\leq7$.
	
	Moreover, by a computer search, we find the sizes for the
	largest maximal partial spreads in $\mathrm{PG}(5, q)$ for $q\leq5$: we prove that the largest maximal partial line
	spreads have deficiency  $\delta=2$  for $q=2$, $\delta=\frac{q+3}{2}$ for $q=3, 5$ and $\delta=\sqrt{q}+1$ for
	$q=4$.
	
	Again by a computer search, we find many new results for $q\leq5$.
	
\section{Known results}
	Maximal partial line spreads (from now on MPS) in $\mathrm{PG}(5, q)$ have
	been investigated by several authors.
	
	In \cite{beut1}, Beutelspacher proved the existence of MPS of size $q^3+q^2+1$.  
	In \cite{gova},  Govaerts showed that this size is the smallest size for MPS
	in $\mathrm{PG}(5, q)$.

	We recall that a spread in 	$\mathrm{PG}(N, q)$, $N$ odd, have size
	\[
	\frac{q^{N+1}-1}{q^{2}-1}.
	\]
	
	Also, we  have the following results.
	
	\begin{theorem}\label{teo2}
		\emph{(\cite{gaszon1})} For any
		$q$ odd and for $q$ even, $q>q_{0}$, in $\mathrm{PG}(N, q),$ $N\geq5$
		odd, there exist maximal partial spreads of any size between $9Nq^{N-2}\log q$
		and $\frac{q^{N+1}-1}{q^2-1}-q+1$.
	\end{theorem}
	
	We remark that the above interval is empty for $q\leq 243$. So, for $q\leq 7$, the only
	known sizes are $q^3+q^2+1$ and the size of a spread.
	
	Upper bounds
	for the size of the largest example of MPS
	were achieved by Govaerts and Storme \cite{govastor2,govastor1}.

	If $\mathcal{S}$ is a partial line spread
	of $\mathrm{PG}(N, q)$, $N$ odd, and $\left|\mathcal{S}\right|=\frac{q^{N+1}-1}{q^{2}-1}-\delta$, 
	then we say that $\mathcal{S}$ has \emph{deficiency} $\delta$. 
		\begin{theorem}\label{epsilon} \emph{(\cite{govastor1})}
	Let $\epsilon=2$ in the case $q=2$ and let $q+\epsilon$ be the
	size of the smallest non-trivial blocking sets in $\mathrm{PG}(2, q)$
	in the case $q>2$. Suppose $N$ odd and $\delta<\epsilon$. Unless $\delta=0$
	there exist no maximal partial line spreads with deficiency $\delta$
	in $\mathrm{PG}(N, q)$.
	\end{theorem}
	
	\begin{corollary} \label{delta}
	\emph{(\cite{govastor1})}	Let $\mathcal{S}$ be a partial $s$-spread of $\mathrm{PG}(N, q)$, $N$ odd,
		of deficiency $\delta>0$. Then
		\begin{enumerate}
			\item \emph{(\cite{bruen1})} $\delta\geq\sqrt{q}+1$ when $q$ is square;
			\item  \emph{(\cite{blok2})}  $\delta\geq c_{p}q^{2/3}+1$, when $q=p^{h}$, $h$ odd, $h>2$, $p$
			prime, $c_{2}=c_{3}=2^{-1/3}$ and $c_{p}=1$ for $p>3$;
			\item \emph{(\cite{blok1})} $\delta\geq\left(q+3\right)/2$ when $q$ is an odd prime.
		\end{enumerate}	
	\end{corollary}

	\begin{theorem}\emph{(\cite{govastor2})}
		Let $\mathcal{S}$ be a maximal partial line spread of  $\mathrm{PG}(N, q)$, $N$ odd,
		 $q>16$, of deficiency $0<\delta<q^{5/8}/\sqrt{2}+1$. Then $\delta\equiv0$
		($\operatorname{mod}\sqrt{q}+1$) and the set of holes of $\mathcal{S}$ is the disjoint
		union of subgeometries $\mathrm{PG}(2s+1, \sqrt{q})$.\newline
		Moreover, $\delta\geq2\left(\sqrt{q}+1\right)$ when $q>4$.
	\end{theorem}

\section{A geometric construction of maximal partial line spreads in $\mathrm{PG}(5, q)$}

We start by proving the following lemma, similarly to Lemma 2.1 in \cite{iura3}. 

\begin{lemma}\label{lemma}
	
	In $\mathrm{PG}(5, q)$, $q$ a prime power, let $S$ be a hyperplane
	and $X$ a point of $S$. 
	Let $\mathcal{L}$ be a set of lines not of $S$, not through $X$, not two of them meeting outside $S$, and such that $\left|\mathcal{L}\right|<q^{3}$.
	Then there is a line through $X$ not of $S$ and skew to every line of $\mathcal{L}$.
	
\end{lemma}

Let $\mathcal{H}$ be a hyperplane of $\mathrm{PG}(5, q)$ and let
$\mathcal{F}$ be a largest partial line spread of $\mathcal{H}$,
the size of which is $q^{3}+1$. So, in $\mathcal{H}$ there are exactly
$q^{2}$ \emph{ holes}, that is points not over the lines of $\mathcal{F}$.
By Lemma \ref{lemma} it follows that there exists a set of $q^{3}$ mutually
disjoint lines of $\mathrm{PG}(5, q)$, not of $\mathcal{H}$, covering
any set of $q^{3}$ points of $\mathcal{H}$. So, we choose a line
set $\mathcal{F}'$ of $q^{2}$ lines, not in $\mathcal{H}$, covering the  $q^{2}$
holes of $\mathcal{H}$. In this way we get a maximal partial spread $\mathcal{F}\cup\mathcal{F}'$
of $\mathrm{PG}(5, q)$, with size $q^{3}+q^{2}+1$. At this point
we deprive $\mathcal{F}$ of a line $r_{1}$ and cover the points
of $r_{1}$ through $q+1$ lines, $r_{1}^{1}, r_{1}^{2},\ldots, r_{1}^{q+1}$,
not of $\mathcal{H}$, whose existence is guaranteed by above lemma,
obtaining the MPS
\[
\mathcal{F}_{1}=\left(\left(\mathcal{F}\cup\mathcal{F}'\right)-r_{1}\right)\cup\left\{ r_{1}^{1}, r_{1}^{2},\ldots, r_{1}^{q+1}\right\} .
\]
We proceed depriving $\mathcal{F}_1$ of a line $r_{2}$ of $\mathcal{F}$ and,
by using the above lemma again, we find $q+1$ mutually skew lines
$r_{2}^{1}, r_{2}^{2},\ldots, r_{2}^{q+1}$, not of $\mathcal{H}$ and
covering $r_{2}$. The line set 
\[
\mathcal{F}_{2}=\left(\mathcal{F}_{1}-r_{2}\right)\cup\left\{ r_{2}^{1}, r_{2}^{2},\ldots, r_{2}^{q+1}\right\} 
\]
is a MPS of size $q^{3}+q^{2}+2q+1$.

We proceed up to deprive $\mathcal{F}$ of $n=\left\lfloor \frac{q^{3}-q^{2}}{q+1}\right\rfloor$
lines, obtaining a MPS  
\[
\mathcal{F}_{n}=\left(\mathcal{F}_{n-1}-r_{n}\right)\cup\left\{ r_{n}^{1}, r_{n}^{2},\ldots, r_{n}^{q+1}\right\},
\]
with $q^{3}+q^{2}+nq+1$ lines.

So, we give a class of maximal partial line spreads of sizes $q^{3}+q^{2}+kq+1$,
with $1\leq k\leq\left\lfloor \frac{q^{3}-q^{2}}{q+1}\right\rfloor$.

So we get the following theorem.

\begin{theorem}\label{teo2.2}

In $\textnormal{PG}(5, q)$, $q$ a prime power, there are 
maximal partial line spreads of size  $q^{3}+q^2+kq+1$, for every integer $ k=1, 2,\ldots, \left\lfloor \frac{q^{3}-q^{2}}{q+1}\right\rfloor$.
\end{theorem}

We note that, for $N=5$ and $q\geq2$, 

\[
9Nq^{N-2}\log q>q^{3}+q^2+\left\lfloor \frac{q^{3}-q^{2}}{q+1}\right\rfloor q+1,
\] 
and so the theorem gives 
$\left\lfloor \frac{q^{3}-q^{2}}{q+1}\right\rfloor$ unknown cardinalities for any $q$.

\section{Computer search of  maximal partial line spreads in $\mathrm{PG}(5, q)$}

We firstly report, in the Table \ref{tab1}, for the investigated values of $q$, the minimum sizes $q^{3}+q^{2}+1$, the maximum sizes (see Theorem \ref{epsilon} and Corollary \ref{delta}), the size $\frac{q^{6}-1}{q^{2}-1}-q+1$ (see Theorem \ref{teo2}) and the sizes of the spreads $\frac{q^{6}-1}{q^{2}-1}$.

	\begin{table}[h]\begin{center}
		\caption{} \label{tab1}
\begin{tabular}{cccccc}
	\toprule
	$q$ & $q^{3}+q^{2}+1$ & $\delta$  & Max size & $\frac{q^{6}-1}{q^{2}-1}-q+1$ & Spread\tabularnewline
	\hline 
	$2$ & $13$ & $\geq2$ & $\leq19$ & $20$ & $21$\tabularnewline
	\hline 
	$3$ & $37$ & $\geq3$ & $\leq88$ & $89$ & $91$\tabularnewline
	\hline 
	$4$ & $81$ & $\geq3$ & $\leq270$ & $270$ & $273$\tabularnewline
	\hline 
	$5$ & $151$ & $\geq4$ & $\leq647$ & $647$ & $651$\tabularnewline
	\hline 
	$7$ & $393$ & $\geq5$ & $\leq2446$ & $2445$ & $2451$\tabularnewline
	\bottomrule
\end{tabular}\end{center}
\end{table}

In the Table \ref{tab2}, we report the obtained cardinalities of type $q^{3}+q^{2}+k\cdot q+1$ of MPS, for $q\leq 7$, by the Theorem \ref{teo2.2} and by the first algorithm for the computer search.

In the Table \ref{tab2}, $k_{min}=1$ and
$\left|\mathcal{F}_{min}\right|=q^{3}+q^{2}+q+1$; $k_{max}$ is the maximum value of $k$ and $\left|\mathcal{F}_{max}\right|=q^{3}+q^{2}+k_{max}\cdot q+ 1$.

\begin{table}[h]\begin{center}
	\caption{} \label{tab2}
\begin{tabular}{cccc}
	\toprule
	\multirow{2}{*}{$q$} & \multicolumn{2}{c}{Theorem \ref{teo2.2}} & \multicolumn{1}{c}{Computer search}
	\tabularnewline
	& $\left|\mathcal{F}_{min}\right|$ & $k_{max}\rightarrow\left|\mathcal{F}_{max}\right|$ & $k_{max}\rightarrow\left|\mathcal{F}_{max}\right|$\tabularnewline
	\hline 
	$2$ & $15$ & $1\rightarrow15$ & $3\rightarrow19$\tabularnewline
	\hline 
	$3$ & $40$ & $4\rightarrow49$ & $12\rightarrow73$\tabularnewline
	\hline 
	$4$ & $85$ & $9\rightarrow117$ & $28\rightarrow193$\tabularnewline
	\hline 
	$5$ & $156$ & $16\rightarrow231$ & $59\rightarrow446$\tabularnewline
	\hline 
	$7$ & $400$ & $36\rightarrow645$ & $163\rightarrow1534$\tabularnewline
	\bottomrule
\end{tabular}\end{center}
\end{table}

So we get the following theorem:

\begin{theorem}In $\mathrm{PG}(5, q)$, $q\leq7$, there are maximal partial line
	spreads of size $q^{3}+q^{2}+kq+1$, for every integer $k=1, 2,\ldots,  
	\left\lfloor\frac{4}{9}q^3\right\rfloor$.
\end{theorem}

We now report all the obtained cardinalities by computer search, that we denote by
$b+k_{m}^{n}\cdot q$, where $b=q^3+q^2+1$ and $m\leq k\leq n$, $k$ integer.
The underlined
values are those not known and that cannot be obtained through the Theorem \ref{teo2.2}.

\begin{description}
	\item [- $q=2$]:  $13,$ $15$, $\underline{16-19}$, $21$.
	\item  [- $q=3$]: $37+k_{0}^{4}\cdot3$, $\underline{37+k_{5}^{7}\cdot3}$,
	 $\underline{60-88}$, $91$.
	\item  [- $q=4$]: $81+k_{0}^{9}\cdot4$, $\underline{81+k_{10}^{19}\cdot4}$,
	$\underline{159-270}$, $273$.
	\item  [- $q=5$]: $151+k_{0}^{16}\cdot5$, $\underline{151+k_{17}^{49}\cdot5}$,
 $\underline{399-647}$,  $651$.
	\item [- $q=7$]: $393+k_{0}^{36}\cdot7$, $\underline{393+k_{37}^{163}\cdot7}$.
\end{description}

We note that all smallest MPS have sizes of type $q^3+q^2+kq+1$, while, starting from
a certain size, there are all the possible sizes.

By the previous result, we get the following theorem (see also Table \ref{tab1}).
\begin{theorem}
In $\mathrm{PG}(5, q)$, $q\leq5$, the largest maximal partial line
spreads have deficiency  $\delta=2$  for $q=2$, $\delta=\frac{q+3}{2}$ for $q=3, 5$ and $\delta=\sqrt{q}+1$ for
$q=4$.
\end{theorem}
The Theorem \ref{teo2}, for $N=5$ and $q=3, 4, 5$,  can be modified as follows.
	\begin{theorem}
 In $\mathrm{PG}(5, q)$, 
 $q=3, 4, 5$, there exist maximal partial line spreads of any size between $2q^{3}\log q$
		and $\frac{q^{6}-1}{q^2-1}-\delta$.
	\end{theorem}

\begin{center}
--------------------
\end{center}
\noindent Maurizio Iurlo\\
\noindent Largo dell'Olgiata, 15\\
\noindent 00123 Roma \\
\noindent Italy\\
\noindent maurizio.iurlo@istruzione.it\\
\noindent http://www.maurizioiurlo.com \\

\end{document}